\documentclass[14pt,reqno,intlim]{amsart}

\usepackage{amssymb,amsmath}

\usepackage[arrow,matrix,curve]{xy}

\newdir{ >}{{}*!/-5pt/\dir{>}}

\newcommand{\nc}{\newcommand}

\nc{\bC}{\bold{C}} \nc{\bN}{\Bbb{N}} \nc{\cF}{\mathcal{F}}
\nc{\cE}{\mathcal{E}} \nc{\cR}{\mathcal{R}} \nc{\cM}{\mathcal{M}}
\nc{\al}{\alpha} \nc{\bt}{\beta} \nc{\gm}{\gamma} \nc{\dl}{\delta}
\nc{\om}{\omega} \nc{\sg}{\sigma} \nc{\Sg}{\Sigma} \nc{\vf}{\varphi}
\nc{\ve}{\varepsilon} \nc{\os}{\overset} \nc{\ol}{\overline}
\nc{\ul}{\underline} \nc{\us}{\underset} \nc{\sbs}{\subset}
\nc{\bsl}{\backslash} \nc{\Ra}{\Rightarrow}
\nc{\lra}{\longrightarrow} \nc{\all}{\allowdisplaybreaks}

\nc{\Codes}{\operatorname{{\bold{Codes}}}}
\nc{\RegMono}{\operatorname{\mathcal{R}{\rm{eg}\mathcal{M}{\rm{ono}\!}}}}
\nc{\RegEpi}{\operatorname{\mathcal{R}{\rm{eg}\mathcal{E}{\rm{pi}\!}}}}
\nc{\Mn}{\operatorname{\mathcal{M}{\rm{ono}\!}}}
\nc{\Ep}{\operatorname{\mathcal{E}{\rm{pi}\!}}}
\nc{\Rg}{\operatorname{\mathcal{R}{\rm{eg}\!}}}
\nc{\Ob}{\operatorname{Ob\!}}

\numberwithin{equation}{section}

\newtheorem{theo}{\ \ \ Theorem}[section]

\newtheorem{prop}[theo]{\ \ \ Proposition}
\newtheorem{cor}[theo]{\ \ \ Corollary}

\theoremstyle{definition}

\theoremstyle{remark}

\begin{document}

\title[]
{Hereditary QF-$3^{+}$ rings}

\author{Dali Zangurashvili}

\maketitle

\begin{abstract}
With the aid of the torsion-theoretical approach, several properties/characterizations of hereditary QF-$3^{+}$ rings are found. These characterizations are formulated in terms of the category of projective modules, the category of injective projective modules, and also in terms of the maximal left/right ring of quotients by Utumi. Moreover, it is shown that, for a hereditary QF-$3^{+}$ ring, the `largest stable submodule' radical (on the category of left/right modules) is permutable with injective envelopes. Besides, it is shown that there is a bimorphism (in the category of associative rings with identity) from such a ring to a semisimple left/right Artinian ring. 
\vskip+2mm
\noindent{\bf Key words and phrases}:  QF-$3^{+}$ ring, hereditary ring, torsion theory, projective module, injective envelope, stable module.
\vskip+2mm

\end{abstract}

\section{Introduction}


The notion of a left/right QF-$3^{+}$ ring was introduced by Colby and Rutter \cite{CR1} as a generalization of that of a left/right  QF-3 algebra by Thrall \cite{T}. A left (resp. right) QF-$3^{+}$ ring is defined as a ring $\Lambda$ the injective envelope of which, viewed as a left (resp. right) module $_{\Lambda}\Lambda$ (resp. $\Lambda_{\Lambda}$) over itself, is projective. Colby and Rutter proved that, if a ring is left hereditary, then the notion of a QF-$3^{+}$ ring is left-right symmetric; moreover, a left hereditary QF-$3^{+}$ ring is also right hereditary \cite[proof of Theorem 3.2, p. 384]{CR}; so we drop the prefix `left/right', where it is appropriate. 

A number of characterizations of hereditary QF-$3^{+}$ rings are known (see, e.g., \cite{WMJ}, \cite{CR},  \cite{MZ}). In \cite{Z}, we proved that a left hereditary ring is QF-$3^{+}$ if and only if the pair of left module classes 
\begin{equation}
(Stable\ modules,\ Projective\ modules) 
\end{equation}
\noindent is a torsion theory (recall that a module is called stable if it has no nonzero projective direct summands. In the case of a left hereditary ring, a left module is stable if and only if it has the zero dual \cite[Lemma 2.6]{MZ}). It is easy to see (Corollary \ref{st7}) that the left modules which are divisible with respect to this torsion theory are precisely the injective modules, and the divisible envelope of any projective module is its injective envelope. As it was noted in \cite[Remark 5.6]{Z}, torsion theory (1.1) coincides with the largest torsion theory for which the left module $_\Lambda \Lambda$ is torsionfree, considered by Lambek in \cite{L2}.

The present paper is devoted to some implications of the latter facts. Namely, new properties and characterizations of hereditary QF-3$^{+}$ rings are given. 
 Several of them are formulated in terms of the category of projective modules, and complement our recent result asserting that a left hereditary ring $\Lambda$ is left perfect and right coherent\footnote{Note that the class of such rings contains the class of hereditary QF-3$^{+}$ rings (see, e.g., \cite[Proposition 9.2]{MZ}).} if and only if this category is a reflective subcategory of the category $\Lambda$-$Mod$ of all left $\Lambda$-modules \cite[Theorem 5.2]{Z}. In the present paper, it is shown that a left hereditary ring is QF-3$^{+}$ if and only if the category of projective left modules is reflective in the category $\Lambda$-$Mod$ and the reflector preserves monomorphisms (Theorem \ref{t1.1}). In view of these results, note also that this reflector is left exact if and only if the ring is semisimple (Proposition \ref{st3}). Moreover, it is shown that, for a hereditary QF-3$^{+}$ ring, the category of stable modules is Abelian, while the category of projective modules is not such unless the ring is semisimple.

Further, it is shown that, for a hereditary QF-3$^{+}$ ring, the injective envelope of a projective left module can be turned into a functor (on the category of such modules), so that the embedding of a projective left module into its injective envelope becomes natural (Proposition \ref{l1}). It is shown that this functor is permutable with the `largest stable submodule' radical (Corollary \ref{st5}). 

Moreover, employing the same torsion-theoretic approach, some further properties of hereditary QF-3$^{+}$ rings and their maximal left rings of quotients $Q^{l}_{max}(\Lambda)$ (in the sense of \cite{U} by Utumi) are obtained. In particular, a new proof of the well-known result asserting that a ring structure that extends the left $\Lambda$-module structure can be introduced on the injective envelope of a hereditary QF-$3^{+}$ ring, and the ring obtained in this way is isomorphic to $Q^{l}_{max}(\Lambda)$ (see, e.g., \cite[the proof of Theorem 3,2]{CR}). With the aid of this fact, it is shown that  there is a bimorphism (in the category of associative rings with identity), i.e., a morphism that is both monomorphism and epimorphism, from such a ring to a semisimple left Artinian ring (Corollary 4.10). 

Further, criteria for a left hereditary ring to be QF-$3^{+}$ are given in terms of their maximal left rings of quotients $Q^{l}_{max}(\Lambda)$. One of them asserts that a left hereditary ring $\Lambda$ is QF-$3^{+}$ if and only if $Q^{l}_{max}(\Lambda)$ is semisimple and projective as a left $\Lambda$-module (Theorem 4.11; observe that the latter condition can be obtained from the definition of a left QF-$3^{+}$ ring by replacing the injective envelope of a ring by its left maximal ring of quotients). Another criterion states that a left hereditary ring $\Lambda$ is QF-$3^{+}$ if and only if the \textit{right} $\Lambda$-module $Q^{l}_{max}(\Lambda)_{\Lambda}$ is projective and any left $Q^{l}_{max}(\Lambda)$-module (including the one $_{\Lambda}Q^{l}_{max}(\Lambda)$) is projective as a left $\Lambda$-module (Theorem 4.12). 

Finally, note that the `right' versions of all above-mentioned statements are valid as well.

Financial support from  Shota Rustaveli  National Science Foundation of Georgia
(Ref.: FR-24-8249) is gratefully acknowledged.

\section{Preliminaries}

Throughout the paper, $\Lambda$ is an associative ring with identity, and `module' means a left unital module over $\Lambda$. Let $\Lambda$-$Mod$ be the category of such modules. When no confusion might arise, we use the same symbol for a class of modules and the corresponding full subcategory of the category $\Lambda$-$Mod$. Namely, the symbols $\Lambda$-$Mod_{St}$, $\Lambda$-$Mod_{Pr}$, and $\Lambda$-$Mod_{PrInj}$ denote the classes/categories of resp. stable, projective, and projective injective modules.

Moreover, as usual, the symbol $_{\Lambda}\Lambda$ denotes the ring $\Lambda$, viewed as a left module over itself.

For a module $M$, the symbol $I(M)$ denotes the injective envelope of a module $M$, while $\textbf{R}(M)$ denotes the largest stable submodule of $M$ (provided that it exists). In the case of a left hereditary ring, its existence follows from \cite[Theorem 1]{J} and \cite[Proposition 2.6]{MZ}. In that case, $\textbf{R}$ is a radical, as it is pointed out in \cite[Proposition 3.3]{Z}. We use the notion of a radical as well as those of a pre-torsion theory and torsion theory in the sense of the book \cite{L2} by Lambek. We assume that a reader is familiar with the basics of torsion theory presented in this book. Below we recall only few notions and facts from this theory.

Let $(\mathcal{B}, \mathcal{C})$ be a torsion theory, and $\mathbf{T}$ be the corresponding radical. A module $M$ is called divisible if $I(M)/M$ is torsionfree, i.e., lies in $\mathcal{C}$. The divisible envelope of a module is defined as a unique submodule  $D(M)$ of $I(M)$ such that $\mathbf{T}(I(M)/M)=D(M)/M$.

\begin{prop} \cite[Proposition 0.7]{L2} \label{p1}
For any module $M$, the following conditions hold:
 
 (1) $M\rightarrowtail D(M)$ is an essential extension; 
 
 (2) D(M)/M is torsion, i.e., lies in the class $\mathcal{B}$; 
 
 (3) $D(M)$ is divisible. 
 
\noindent Moreover, the extension  $M\rightarrowtail D(M)$ is detrmined by these conditions uniquely up to isomorphism.
\end{prop}

The mapping $M\mapsto D(M)$ can be turned into a functor $\mathbf{D}$ from the category of torsionfree modules to the category of divisible torsionfree modules such that the embedding $M\rightarrowtail \mathbf{D}(M)$ is natural; this can be done in a unique way.

\begin{theo} \cite[Proposition 0.8]{L2} \label{2.2}
The category $\mathcal{A}$ of divisible torsionfree modules is Abelian and reflective in the category of torsionfree modules (with the reflector $\mathbf{D}$) and also in $\Lambda$-$Mod$ with the reflector $\mathbf{Q}$ given by the formula:
$$\mathbf{Q}(M)=D(M/\mathbf{T}(M)).$$
The functor $\mathbf{Q}$ is exact. 
\end{theo}

The left module $\mathbf{Q}(_\Lambda \Lambda)$ can be equipped with a structure of an associative ring with identity due to the group isomorphisms 
\begin{equation} \label{2.1}
\mathbf{Q}(_\Lambda \Lambda)\simeq Hom_{\Lambda}(_\Lambda\Lambda, \mathbf{Q}(_\Lambda \Lambda))\simeq Hom_{\Lambda}(\mathbf{Q}(_\Lambda \Lambda), \mathbf{Q}(_\Lambda \Lambda)).
\end{equation}
This ring is called the \textit{left ring of quotients} of $_\Lambda \Lambda$ with respect to the torsion theory  $(\mathcal{B}, \mathcal{C})$ \cite[Example 5, page 25]{L2}. 

Throughout the paper, we will repeatedly employ the following statement.

\begin{theo} \cite[Theorem 5.2, Remark 5.4]{Z} \label{z1}
For a ring $\Lambda$, the  conditions (1)-(8) below are equivalent. When these conditions are satisfied, for the submodule $P$ mentioned in the condition (4), there is an isomorphism $i:P\rightarrow \mathbf{r}(M)$, where $\mathbf{r}$ is the reflector 
\begin{center}
$\Lambda$-$Mod\rightarrow \Lambda$-$Mod_{Pr}$
\end{center}
that exists due to the condition (2). Moreover, the unit of this reflection is the composition of the canonical projection $\mathbf{R}(M)\oplus P\twoheadrightarrow P$ with the isomorphism $i$.
\begin{enumerate}
\item The ring $\Lambda$ is left hereditary, left perfect and right coherent;\vskip+2mm

\item the subcategory $\Lambda$-$Mod_{Pr}$ of $\Lambda$-$Mod$ is epireflective;
\vskip+2mm 

\item the ring $\Lambda$ is left hereditary and the subcategory $\Lambda$-$Mod_{Pr}$ of $\Lambda$-$Mod$ is reflective;\vskip+2mm

\item the ring $\Lambda$ is left hereditary and, moreover, for any module $M$, there exists a representation $M= \textbf{R}(M)\oplus P$ with a projective submodule $P$;\vskip+2mm

\item \textit{any module can be represented as a direct sum of a stable submodule and a projective submodule;}\vskip+2mm

\item \textit{any module can be represented as a direct sum of a stable submodule and a projective submodule, and such a representation is unique up to isomorphism (in the sense that if $M= S\oplus P$ and $M=S'\oplus P'$ with  stable submodules $S$ and $S'$ and projective submodules $P$ and $P'$, then $S$ is equal to $S'$, and $P$ is isomorphic to $P'$)}; \vskip+2mm

\item \textit{the pair of module classes $(\Lambda$-$Mod_{St},\Lambda$-$Mod_{Pr})$ is a pre-torsion theory;}\vskip+2mm

\item \textit{the ring $\Lambda$ is left hereditary and, moreover, for any module $M$, $\mathbf{R}(M)=0$ if and only if $M$ is projective}. \vskip+2mm

\end{enumerate}
\end{theo}

Finally, recall the following statement.
\begin{theo} \cite[Proposition 9.2, p. 58]{K}\label{k2}
Let $(\mathcal{B}, \mathcal{C})$ be a torsion theory, and $T$ be the corresponding radical. The class $\mathcal{B}$ is closed under essential extensions if and only if $$\mathbf{T}(I(M))=I(\mathbf{T}(M)),$$ for any module $M$.
\end{theo}

\section{The categories of stable and projective modules over a hereditary QF-$3^{+}$ ring}
Recall the following statement.

\begin{theo} \cite[Theorem 5.6]{Z} \label{t1.1}
A ring $\Lambda$ is a hereditary  QF-3$^+$ ring if and only if the pair of module classes 
\begin{equation}
(\Lambda\text{-}Mod_{St}, \Lambda\text{-}Mod_{Pr})
\end{equation}

\noindent is a torsion theory.

\end{theo}

Theorem \ref{t1.1} immediately implies the following statement.
\begin{theo} \cite[Theorem 2.3]{CR}\label{t1.1.1}
Let $\Lambda$ be a left hereditary ring. It is a QF-3$^+$ ring if and only if the class $Mod_{Pr}$ is closed under injective envelopes. 

\end{theo}

We provide more criteria for a left hereditary ring to be a QF-$3^{+}$ ring.

\begin{theo}\label{t1}
Let $\Lambda$ be a left hereditary ring. The following conditions are equivalent:\vskip+1mm
\begin{enumerate}
\item the ring $\Lambda$ is a QF-$3^{+}$ ring;\vskip+1mm

\item the ring $\Lambda$ is left perfect and right coherent, and, moreover, the category $\Lambda$-$Mod_{St}$ is a Serre subcategory of $\Lambda$-$Mod$;\vskip+1mm

\item the ring $\Lambda$ is left perfect and right coherent, and, moreover, the category $\Lambda$-$Mod_{St}$ is a localizing subcategory of $\Lambda$-$Mod$;\vskip+1mm

\item the ring $\Lambda$ is left perfect and right coherent, and, moreover, for any module $M$, we have 
\begin{equation} \label{3.1}
\textbf{R}(M)=\cap_{f:M\rightarrow I(_\Lambda \Lambda)} Ker f;
\end{equation}

\item the ring $\Lambda$ is left perfect and right coherent, and, moreover, no stable module contains a nonzero projective submodule;\vskip+1mm

\item the subcategory $\Lambda$-$Mod_{Pr}$ of $\Lambda$-$Mod$ is reflective, and, moreover, the reflector $$\textbf{r}:\Lambda\textit{-}Mod\rightarrow \Lambda\textit{-}Mod_{Pr}$$
\noindent preserves monomorphisms;\vskip+1mm

\item the subcategory $\Lambda$-$Mod_{Pr}$ of $\Lambda$-$Mod$ is reflective, and, moreover, the image of the embedding $_\Lambda \Lambda\rightarrowtail I(_\Lambda \Lambda)$ under the reflector $\textbf{r}$ is a monomorphism.
\vskip+1mm
\end{enumerate}
\end{theo}

\begin{proof}

(1)$\Rightarrow$(2): A left hereditary QF-$3^{+}$ ring is semi-primary and left Noetherian, as is shown in \cite[the proof of Theorem 3.2]{CR}. This immediately implies that it is left perfect and right coherent (see also of \cite[ Proposition 9.2(a)]{MZ}). Now it suffices to apply Theorem \ref{t1.1}. 

The equivalence (2)$\Leftrightarrow$(3) follows from \cite[Lemma 3.1(b)]{Z}. 

(3)$\Rightarrow$(1): The class $\Lambda$-$Mod_{St}$ is closed under isomorphic images, factor modules, extensions, direct sums and submodules. Hence the pair $(\Lambda$-$Mod_{St}, (\Lambda$-$Mod_{St})^{\perp})$ is a torsion theory. According to \cite[Theorem 5.2]{Z}, we have $(\Lambda$-$Mod_{St})^{\perp}=\Lambda$-$Mod_{Pr}$. Now the claim follows from Theorem \ref{t1.1}.

(1)$\Rightarrow$(4): Since $\mathbf{R}(M)$ is stable, it is contained in $Ker\; f$, for any $f:M\rightarrow I(_\Lambda\Lambda)$, by \cite[Lemma 3.2]{MZ}. On the other hand, the right-hand part $N$ of (\ref{3.1}) is stable. Indeed, for a homomorphism $f:N\rightarrow_{\Lambda}\Lambda$, there is a homomorphism $g$ such that the square 
\begin{equation}
\begin{xymatrix}
{N\ar@{>->}[r]\ar[d]_{f}&M\ar@{-->}[d]^{g}\\
_\Lambda\Lambda\ar@{ >->}[r]&I(_\Lambda\Lambda)}
\end{xymatrix}
\end{equation}
is commutative. Obviously, $N$ is contained in the kernel of $g$. Therefore, $f=0$. Applying \cite[Lemma 3.2]{Z} again, we obtain that $N$ is stable, and hence is a submodule of $\mathbf{R}(M)$.

(4)$\Rightarrow$(1): We have $\mathbf{R}(I(_\Lambda\Lambda))=0$. Then, according to Theorem \ref{z1}, $I(_\Lambda\Lambda)$ is projective.

The equivalence (1)$\Leftrightarrow$(5) follows from Theorems 2.3 and 3.1 and the fact that any hereditary QF-$3^{+}$ ring is left perfect and right coherent (see the proof of (1)$\Leftrightarrow$(2)).

(1)$\Rightarrow$(6): We already know that $\Lambda$ is left perfect and right coherent. Theorem 2.3 implies that the subcategory $\Lambda$-$Mod_{Pr}$ of $\Lambda$-$Mod$ is epi-reflective. Let $\textbf{r}$ be the corresponding reflector and $\eta$ be the unit of the reflection. Let $f:M\rightarrowtail N$ be a monomorphism in $\Lambda$-$Mod$. Theorem \ref{t1.1} implies that the kernels of $\eta_M$ and $\eta_N$ are resp. $\mathbf{R}(M)$ and $\mathbf{R}(N)$. Consider a homomorphism $g:K\rightarrow \textbf{r}(M)$ with $\mathbf{r}(f)g=0$.  Let $j$ be a section of $\eta_M$. The equality $\eta_Nfj g=0$ implies that there is a morphism $h$ with $\iota_N h=fj g$. Theorem \ref{t1.1} implies that $\mathbf{R}(M)=
M\cap \mathbf{R}(N)$. Therefore, the left-hand square in the following diagram is a pullback:
\begin{equation}
\begin{xymatrix}
{&K\ar[dr]^{g}\ar@{-->}[ddl]^{h}\ar@{-->}[dl]_{d}\\
\textbf{R}(M)\ar@{ >->}[r]^{\iota_M}\ar@{ >->}[d]_{\textbf{R}(f)}&M\ar@{ ->>}[r]^{\eta_M}\ar@{ >->}[d]^{f}&\textbf{r}(M)\ar[d]^{\textbf{r}(f)}\\
\textbf{R}(N)\ar@{ >->}_{\iota_N}[r]&N\ar@{ ->>}[r]_{\eta_N}&\textbf{r}(N)}
\end{xymatrix}
\end{equation}
\noindent Hence there is a morphism $d$ such that $\iota_Md=j g$. This implies that $\eta_M\iota_Md=\eta_Mj g=g$, and hence $g=0$. 

The implication (6)$\Rightarrow$(7) is obvious.

(7)$\Rightarrow$(1): Consider the commutative diagram
\begin{equation}
\begin{xymatrix}
{_\Lambda\Lambda\ar[r]^{\simeq}\ar@{ >->}[d]_{i}&\textbf{r}(_\Lambda\Lambda)\ar@{ >->}[d]^{\textbf{r}(i)}\\
I(_\Lambda \Lambda)\ar_{\eta_{I(_\Lambda \Lambda)}}[r]&\textbf{r}(I(_\Lambda \Lambda))}
\end{xymatrix}
\end{equation}
\noindent Since $\Lambda$-$Mod_{Pr}$ is a reflective subcategory, $\textbf{r}(i)$ is a monomorphism in $\Lambda$-$Mod$. Since the monomorphism $i$ is essential, we obtain that $\eta_{I(_\Lambda\Lambda)}$ is a monomorphism, and hence the module $I(_{\Lambda} \Lambda)$ is projective.

\end{proof}

Theorem \ref{t1.1} and Theorem \ref{t1} imply the following statement.
\begin{cor} \label{st4}
Let $\Lambda$ be a hereditary QF-$3^{+}$ ring. The category $\Lambda$-$Mod_{St}$ of stable modules is Abelian and closed in $\Lambda$-$Mod$ under essential extensions and hence under injective envelopes.
Moreover, in $\Lambda$-$Mod_{St}$ has injective envelopes.
\end{cor}

\begin{proof}
The fact that $\Lambda$-$Mod_{St}$ is Abelian immediately follows from Theorem \ref{t1}. Let $S$ be a stable module, and $m:S\rightarrowtail M$ be an essential extension. Assume $P$ is a projective submodule of $M$. Theorem \ref{t1.1} implies that the class of stable modules is closed under subobjects. Therefore, the submodule $S\cap P$ is both projective and stable, and hence zero. This implies that $P=0$. 

For the last claim, first note that monomorphisms in $\Lambda$-$Mod_{St}$ are monomorphisms in $\Lambda$-$Mod$ too. Indeed, assume that a homomorphism $f:S\rightarrow S'$ with stable $S$ and $S'$ is a monomorphism in $\Lambda$-$Mod_{St}$. Then the $f$'s kernel $Ker\; f$ (in $\Lambda$-$Mod$) has no nonzero stable submodule. Hence $Ker\;f$ is projective, as it follows from Theorem \ref{z1} and the fact that any hereditary QF-$3^{+}$ ring is left perfect and right coherent mentioned in the proof of Theorem \ref{t1}. Theorem \ref{t1} implies that $Ker\;f=0$, and hence $f$ is a monomorphism in $\Lambda$-$Mod$ too. Therefore, the stable modules which are injective in $\Lambda$-$Mod$ are injective also in $\Lambda$-$Mod_{St}$.
\end{proof}

Corollary \ref{st4} and Theorem \ref{k2} imply the following corollary.

\begin{cor}\label{st5}
Let $\Lambda$ be a hereditary QF-$3^{+}$ ring. For any module $M$, we have 
\vskip+2mm
(a) $\textbf{R}(I(M))=I(\textbf{R}(M))$;
\vskip+2mm

(b) If $I$ is an injective module, then so is $\textbf{R}(I)$.
\end{cor}

\begin{prop}\label{st3}
Let a ring $\Lambda$ be left hereditary, left perfect and right coherent. The reflector  $$\textbf{r}:\Lambda\textit{-}Mod\rightarrow \Lambda\textit{-}Mod_{Pr}$$ \noindent is left exact if and only if the ring is semisimple. 
\end{prop}

\begin{proof} The "if" part of the claim is obvious. To show "only if" part, first note that the functor $\textbf{r}$ maps short exact sequences to sequences which are exact in $\Lambda$-$Mod$. Indeed, let 
\begin{equation}
\begin{xymatrix}
{0\ar[r]&K\ar@{ >->}[r]^{m}&N\ar@{->>}[r]^{e}&M\ar[r]&0\\}
\end{xymatrix}
\end{equation}
\noindent be a short exact sequence in $\Lambda\textit{-}Mod$. Let $\eta$ be the unit of the reflection. Since $\eta_M$ is an epimorphism, so is $\textbf{r}(e)$ in $\Lambda$-$Mod$. Moreover, $\textbf{r}(m)$ is a kernel of $\textbf{r}(e)$ (in both categories $\Lambda$-$Mod$ and $\Lambda$-$Mod_{Pr}$) since $\mathbf{r}$ is left exact. 

Now consider an arbitrary module $M$, and  an epimorphism $e:P\twoheadrightarrow M$ with a projective $P$. Let $K$ be the kernel of $e$. Then the sequence 
\begin{equation}
\begin{xymatrix}
{0\ar[r]&\textbf{r}(K)\ar@{  >->}[r]&\textbf{r}(P)\ar@{->>}[r]^{\textbf{r}(e)}&\textbf{r}(M)\ar[r]&0\\}
\end{xymatrix}
\end{equation}
\noindent is exact in $\Lambda$-$Mod$. At that, the homomorphisms $\eta_K$ and $\eta_P$ are isomorphisms. This implies that so is $\eta_M$, and hence $M$ is projective.

\end{proof}

The category $\Lambda$-$Mod_{Pr}$ is additive as a full reflective subcategory of $\Lambda$-$Mod$. We obtain the following statement.

\begin{cor}
Let $\Lambda$ be a hereditary QF-$3^{+}$ ring. The category $Mod_{Pr}$ of projective modules is Abelian if and only if the ring is semisimple.
\end{cor}

\begin{proof}
Recall  that the reflector of a full reflective Abelian subcategory $\mathcal{A}$ of the category $\Lambda$-$Mod$ is left exact if and only if it preserves monomorphisms (see, e.g., the remark on page 13 of \cite{L2}). Now it suffices to apply Theorem \ref{t1} and Proposition \ref{st3}.
\end{proof}

\section{The maximal left ring of quotients of $\Lambda$}
Let $\Lambda$ be a hereditary QF-$3^{+}$ ring. As was observed in \cite[Remark 5.7]{Z}, the torsion theory $(\Lambda\text{-}Mod_{St}, \Lambda\text{-}Mod_{Pr})$ coincides with the largest torsion theory for which the module $_\Lambda \Lambda$ is torsionfree, considered by Lambek in \cite{L2}. As is shown in \cite{L2}, the modules that are divisible with respect to the latter theory are precisely the rationally complete modules in the sense of \cite{FL}. We have also the following characterization.
\begin{cor}\label{st7}
Let $\Lambda$ be a hereditary QF-$3^{+}$ ring. The modules which are divisible with respect to the torsion theory $(\Lambda\text{-}Mod_{St}, \Lambda\text{-}Mod_{Pr})$ are precisely the injective modules. The divisible envelope of any module is its injective envelope.
\end{cor}

\begin{proof}
If $M$ is divisible, then $I(M)/M$ is projective. This implies that $M$ is a direct summand of $I(M)$, and hence $M$ is injective. The converse implication is obvious. Moreover, by Proposition \ref{p1}, the canonical embedding of a module into its divisible envelope is essential. This implies the second part of the claim.
\end{proof}

\begin{cor}
Let $\Lambda$ be a hereditary left QF-$3^{+}$ ring. Then a module is rationally complete if and only if it is injective.
\end{cor}

Let $\Lambda$-$Mod_{PrInj}$ be the full subcategory of projective injective modules of $\Lambda$-$Mod$. Theorem \ref{2.2} and Corollary \ref{st7} imply the following statement.
\begin{cor}\label{st8}
If $\Lambda$ is a hereditary QF-$3^{+}$ ring, then the subcategory $\Lambda$-$Mod_{PrInj}$ is  reflective in the category $\Lambda$-$Mod_{Pr}$ and also in $\Lambda$-$Mod$; the reflector
$$\textbf{r}':\Lambda\textit{-}Mod\rightarrow \Lambda\textit{-}Mod_{PrInj}$$  in the latter reflection is exact. Moreover, the category  $\Lambda$-$Mod_{PrInj}$ is Abelian.
\end{cor}

The reflector  $$\Lambda\textit{-}Mod_{Pr}\rightarrow \Lambda\textit{-}Mod_{PrInj}$$ mentioned in Corollary \ref{st8} sends a projective module to its injective envelope, and the components of the unit of this reflection are essential monomorphisms $P\rightarrowtail I(P)$. This, in particular, implies the following statement.

\begin{prop} \label{l1}
Let $\Lambda$ be a hereditary QF-$3^{+}$ ring. Then, for any projective modules $P$ and $Q$ and any homomorphism $f:P\rightarrow Q$, there is a unique homomorphism $f'$ making the following diagram commutative:
\begin{equation}
\begin{xymatrix}
{P\ar@{>->}[r]^{i_P}\ar[d]_{f}&I(P)\ar@{-->}[d]^{f'}\\
Q\ar@{>->}[r]^{i_Q}&I(Q)}
\end{xymatrix}
\end{equation}
\end{prop}

Thus, there is precisely one way to turn the mapping $I(-)$ into a functor $$\mathbf{I}:\Lambda\textit{-}Mod_{Pr}\rightarrow \Lambda\textit{-}Mod_{PrInj}$$ for which the embedding $P\rightarrowtail \mathbf{I}(P)$ is natural. This functor is precisely the above-mentioned reflector. It is obviously additive.

As it was mentioned in the Preliminaries, for any torsion theory, the module $\mathbf{Q}(_{\Lambda}\Lambda)$ has a ring structure. Applying this fact, we obtain the well-known statement that, for a hereditary QF-$3^{+}$ ring, there is a ring structure on the injective envelope of $_\Lambda\Lambda$ \cite[page 383]{CR}.  Further, it is well-known \cite{L1} that the left ring of quotients of the largest torsion theory for which the module $_\Lambda \Lambda$ is torsionfree is isomorphic to the maximal left ring of quotients $Q^{l}_{max}(\Lambda)$ of $\Lambda$ introduced by Utumi in \cite{U}. As it was mentioned above, this torsion theory coincides with the torsion theory $(\Lambda\text{-}Mod_{St}, \Lambda\text{-}Mod_{Pr})$ that we are dealing with in this paper provided that $\Lambda$ is a hereditary $QF$-$3^{+}$ ring. Therefore, we obtain the following statement. 

\begin{prop} Let $\Lambda$ be a hereditary QF-$3^{+}$ ring. The left ring of quotients of the torsion theory $(\Lambda\text{-}Mod_{St}, \Lambda\text{-}Mod_{Pr})$ is the injective envelope $I(_\Lambda \Lambda)$ of $_\Lambda \Lambda$ equipped with the following multiplication:
 for $x$ and $y$ in $I(_\Lambda\Lambda)$, $$xy=\mathbf{I}(f_y)(x),$$ where $f_y$ is a left-module homomorphism $_\Lambda\Lambda\rightarrow I(_\Lambda\Lambda)$ that sends $r$ to $ry$. The identity of this ring is $i_{_\Lambda\Lambda}(1)$, where $i_{_\Lambda\Lambda}$ is the essential monomorphism $_\Lambda\Lambda\rightarrowtail I(_\Lambda\Lambda)$. The mapping $i_{_\Lambda\Lambda}$ is a ring homomorphism. The constructed ring structure on the injective envelope $I(_\Lambda \Lambda)$ of $_\Lambda \Lambda$ is the unique ring structure that extends $\Lambda$-left-module structure of $I(_\Lambda \Lambda)$.
 
\begin{proof} It is easy to show that the introduced structure is indeed a ring structure that extends $\Lambda$-left-module structure of $I(_\Lambda \Lambda)$, and then to apply \cite[Theorem 13.11(3)]{L}. The last claim also follows from this theorem.
\end{proof}
 \end{prop}
 
  We denote the constructed ring by the symbol $\widetilde{I(_\Lambda\Lambda)}$. 
  

 \begin{prop} \label{10} \label{4.6}
Let $\Lambda$ be a hereditary QF-$3^{+}$ ring. The rings $\widetilde{I(_\Lambda\Lambda)}$ and $Q^{l}_{max}(\Lambda)$ are isomorphic.
\end{prop}




 \cite[Proposition 5]{L1} implies the following statement.

\begin{cor}
Let $\Lambda$ be a hereditary QF-$3^{+}$ ring. The ring $\widetilde{I(_\Lambda\Lambda)}$ is left self-injective.
\end{cor}

Applying Proposition \ref{10}, we obtain the following statement.
\begin{cor}\label{st11}
Let $\Lambda$ be a hereditary QF-$3^{+}$ ring. Then there is an equivalence $\textbf{F}$  
between the category $\Lambda$-$Mod_{PrInj}$ of projective injective $\Lambda$-modules and the category of left modules over the ring $Q^{l}_{max}(\Lambda)$ such that the composition of the reflector $$\textbf{r}': \Lambda\text{-}Mod\rightarrow \Lambda\text{-}Mod_{PrInj}$$ \noindent with $\textbf{F} $ is naturally isomorphic to the functor $$Q^{l}_{max}(\Lambda)\otimes_{\Lambda} (-):\Lambda\text{-}Mod\rightarrow Q^{l}_{max}(\Lambda)\text{-}Mod,$$ \noindent and, moreover, ${(Mod\; i_{_\Lambda\Lambda})}\textbf{F}=\textbf{i}$:
\begin{equation}
\begin{xymatrix}
{\Lambda\text{-}Mod_{PrInj}\ar@{-->}[d]_{\textbf{F}}\ar[dr]^{\textbf{i}}&\\
Q^{l}_{max}(\Lambda)\text{-}Mod\ar[r]^-{Mod\;i_{_\Lambda\Lambda}}& \Lambda\text{-}Mod}
\end{xymatrix}
\end{equation}
\noindent Here $Mod\; i_{_\Lambda\Lambda}$ is the restriction of scalars functor induced by the ring embedding ${i_{_\Lambda\Lambda}}:\Lambda\rightarrowtail Q^{l}_{max}(\Lambda)$, while $\textbf{i}$ is the embedding functor.  
\noindent 
\end{cor}

\begin{proof}
According to \cite[Proposition 1.2]{L2} and Corollary \ref{st8}, the functor $\mathbf{i}$ is representable with $I(_\Lambda\Lambda)$. Now the fact that there is a functor $\mathbf{F}$ such that $Mod\;i_{_\Lambda\Lambda}\textbf{F}=\textbf{i}$ immediately follows from Proposition \ref{4.6} and \cite[Proposition 1.1(5)]{L2}. To show that $\mathbf{F}$ is an equivalence, first recall that  the ring $\Lambda$ is left Noetherian, as it is noticed in the proof of \cite[Theorem 3.2]{CR}. This together with Papp-Bass's criterion for a ring to be left Noetherian (see, e.g., \cite{GW}) implies that the category $\Lambda\text{-}Mod_{PrInj}$ is closed under coproducts, and the functor $\mathbf{i}$ preserves them. Moreover, if $N$ is a submodule of a module $M$ and both modules are projective and injective, then so is the quotient module $M/N$. Hence the category $\Lambda\text{-}Mod_{PrInj}$ is closed under coequalizers, and the functor $\mathbf{i}$ preserves them. Thus, the category $\Lambda\text{-}Mod_{PrInj}$ is closed under colimits, and the functor $\mathbf{i}$ preserves them. Now it suffices to apply \cite[ Proposition 1.3]{L2}.  
\end{proof}

\begin{cor}\label{st12}
Let $\Lambda$ be a hereditary QF-$3^{+}$ ring. Then 

(a) for any left $Q^{l}_{max}(\Lambda)$-module $M$, the canonical epimorphism $$Q^{l}_{max}(\Lambda)\otimes_{\Lambda} M\twoheadrightarrow\; _{\Lambda}M$$ that sends $s\otimes m$ to $sm$  is an isomorphism;

(b) the ring monomorphism $i_{_\Lambda\Lambda}:\Lambda\rightarrowtail \widetilde{I(_\Lambda\Lambda)}$ is an epimorphism in the category of associative rings with identity;

(c) the left module $_{\Lambda}Q^{l}_{max}(\Lambda)$ is projective.

\end{cor}
\begin{proof}
 The claim 
(a) follows from the fact that the functor $Mod\; i_{_\Lambda\Lambda}$ is full and faithful. The claim (b) follows from the well-known criterion for the restriction-of-scalars functor to be full (see, e.g., \cite[Lemma 3]{L1}). The claim (c) follows from Proposition \ref{10}. 
\end{proof}

Corollary \ref{st11}  implies the well-known fact that, for a hereditary QF-$3^{+}$ ring, the maximal left ring of quotients $Q^{l}_{max}(\Lambda)$ is semisimple (since each object of the category $\Lambda$-$Mod_{PrInj}$ is injective (in it)) \cite[the proof of Theorem 3.2, p. 383]{CR}. In fact, the ring $Q^
{l}_{max}(\Lambda)$ is also left Artinian \cite[page 383]{CR}. Therefore, Corollary \ref{st12}(b) implies the following statement.

\begin{cor}\label{st14}
For any hereditary QF-$3^{+}$ ring $\Lambda$, there is a bimorphism (i.e., a morphism that is both an epimorphism and a monomorphism) in the category of associative rings with identity from $\Lambda$ to a semisimple left Artinian ring.
\end{cor}

\begin{theo}
Let $\Lambda$ be a left hereditary ring. It is a QF-3$^{+}$ ring if and only if its maximal left ring of  quotients $Q^l_{max}(\Lambda)$ is semisimple and projective as a left $\Lambda$-module.
\end{theo}

\begin{proof}
"Only if" part follows from Corollary \ref{st12} (c). For "if" part, note that $\Lambda$ is left nonsingular, by Gabriel's Theorem (see, e.g., \cite[Theorem 13.40]{L}). By Johnson's Theorem, $Q^{l}_{max}(\Lambda)=I(\Lambda)$ (see, e.g., \cite[Theorem 13.36]{L}). This implies the claim.
\end{proof}

Finally, we deal with the problem whether the converse of the first claim of Corollary 4.3 is valid. This, in particular, leads to yet another criterion for a left hereditary ring to be a QF-$3^{+}$ ring that is formulated in terms of its maximal left ring of quotients. 

\begin{theo}\label{st13}
Let $\Lambda$ be a left hereditary ring. The following conditions are equivalent:\vskip+2mm
\begin{enumerate}
\item the ring $\Lambda$ is QF-$3^{+}$ ring;\vskip+2mm

\item the full subcategory $\Lambda$-$Mod_{PrInj}$ of injective projective modules is reflective in $\Lambda$-$Mod$, and $\eta_{_\Lambda\Lambda}$ is a monomorphism, where $\eta$ is the unit of the reflection;
\vskip+2mm

\item the full subcategory $\Lambda$-$Mod_{PrInj}$ of injective projective modules is reflective in $\Lambda$-$Mod_{Pr}$, and $\eta_{_\Lambda\Lambda}$ is a monomorphism, where $\eta$ is the unit of the reflection;
\vskip+2mm

\item the category $\Lambda\text{-}Mod_{Pr}$ has a full subcategory $\mathcal{X}$ such that the composition of the embedding functors $$\mathcal{X}\rightarrow \Lambda\textit{-}Mod_{Pr}\rightarrow \Lambda\textit{-}Mod$$ has a left exact left adjoint $\textbf{r}'$ and $\eta_{_\Lambda\Lambda}$ is a monomorphism, where $\eta$ is the unit of the adjunction;
\vskip+2mm

\item the category $\Lambda\text{-}Mod_{Pr}$ has a full subcategory $\mathcal{X}$ such that the composition of the embedding functors 
$$\mathcal{X}\rightarrow \Lambda\textit{-}Mod_{Pr}\rightarrow \Lambda\textit{-}Mod$$
\noindent has a left adjoint $\textbf{r}'$ that preserves monomorphisms and $\eta_{_\Lambda \Lambda}$ is a monomorphism.\vskip+2mm

\item  the right $\Lambda$-module $Q^{l}_{max}(\Lambda)_{\Lambda}$ is projective, and there is a full and faithful functor $$U: Q^{l}_{max}(\Lambda)\textit{-}Mod\rightarrow \Lambda\textit{-}Mod_{Pr}$$ such that $Mod\;i_{_\Lambda\Lambda}=\mathbf{i}'U$, where $\mathbf{i}'$ is the embedding functor $$\Lambda\textit{-}Mod_{Pr}\rightarrow \Lambda\textit{-}Mod;$$

\item the right $\Lambda$-module $Q^{l}_{max}(\Lambda)_{\Lambda}$ is projective, and any left $ Q^{l}_{max}(\Lambda)$-module is projective as a left $\Lambda$-module.
\end{enumerate}
\vskip+2mm
\end{theo}
\begin{proof}
The implications (1)$\Rightarrow$(2), (1)$ \Rightarrow$(3) and (1)$ \Rightarrow$(4) follow from Corollary 4.3. The implication (4)$\Rightarrow$(5) is obvious.

(2)$\Rightarrow$(1) and (3)$\Rightarrow$(1): Since the module $_\Lambda \Lambda$ can be embedded into a projective injective module,  the injective envelope of $_\Lambda \Lambda$ also can be embedded into such a ring, and hence is projective.

(5)$\Rightarrow$(1): Consider the commutative square
\begin{equation}
\begin{xymatrix}
{_{\Lambda}\Lambda\ar@{  >->}[r]^{i}\ar@{ >->}[d]_{\eta_{_\Lambda \Lambda}}&I(_\Lambda \Lambda)\ar@{ ->}[d]^{\eta_{I(_\Lambda \Lambda)}}\\
\textbf{r}'(_\Lambda \Lambda)\ar@{  >->}[r]^{\textbf{r}'(i)}&\textbf{r}'(I(_\Lambda \Lambda))}
\end{xymatrix}
\end{equation}

\noindent Since $i$ is an essential monomorphism, $\eta_{I(_\Lambda \Lambda)}$ is a monomorphism. Therefore the module $I(_\Lambda \Lambda)$ is projective.

The implication (1)$\Rightarrow$(6) follows from Corollary \ref{st11}, 
Proposition \ref{st12}(c), and the fact that the notion of a hereditary $QF$-$3^{+}$ ring is left-right symmetric \cite[the proof of Theorem 3.2, p.384]{CR}. 

(6)$\Rightarrow$(5): Since the right $\Lambda$-module $Q^{l}_{max}(\Lambda)$ is projective, it is flat. Hence the left adjoint of $Mod\;i_{{\Lambda}\Lambda}$ is exact. Let $\eta$ be the unit of the adjunction
$$Q^{l}_{max}(\Lambda)_{\Lambda}\otimes - \dashv Mod\;i_{{\Lambda}\Lambda}.$$
Its component $\eta_{_\Lambda \Lambda}$ is a monomorphism as it follows from the existence of a left $\Lambda$-module monomorphism $_{\Lambda}\Lambda\rightarrowtail_{\Lambda}Q^{l}_{max}(\Lambda)$.

The equivalence (6)$\Leftrightarrow$(7) is obvious.

\end{proof}
\vskip+2mm

\vskip+2mm
\textit{Author's address:}

\noindent \textit{Dali Zangurashvili, A. Razmadze Mathematical Institute of Tbilisi State University},
\textit{Alexidze Str., Lane II, Tbilisi 0193, Georgia, e-mail: dali.zangurashvili@tsu.ge}

\end{document}